\documentclass[english,preprint]{elsarticle}
\usepackage[T2A,T2A,T2A]{fontenc}
\usepackage[utf8]{inputenc}
\usepackage{babel}
\usepackage{amsmath}
\usepackage{amssymb}
\usepackage{graphicx}
\usepackage[unicode=true,
 bookmarks=false,
 breaklinks=false,pdfborder={0 0 1},backref=section,colorlinks=false]
 {hyperref}

\makeatletter

\usepackage{tikz}
\usetikzlibrary{calc}

\usepackage{babel}
\AtBeginDocument{\DeclareFontEncoding{T2A}{}{}}

\usepackage{lineno}

\usepackage{comment}

\modulolinenumbers[10]
\modulolinenumbers[1]

\makeatother

\begin{document}

\title{Local sensitivity of spatiotemporal structures}

\author{I.A. Shepelev }

\ead{igor\_sar@li.ru }

\address{Department of Physics, Saratov State University,83 Astrakhanskaya
Street, Saratov, 410012, Russia}

\author{A.V. Bukh}

\ead{buh.andrey@yandex.ru }

\address{Department of Physics, Saratov State University,83 Astrakhanskaya
Street, Saratov, 410012, Russia}

\author{S. Ruschel}

\ead{ruschel@math.tu-berlin.de}

\address{Institute of Mathematics, Technische Universität Berlin, Straße des
17. Juni 136, 10623, Berlin, Germany.}

\author{S. Yanchuk}

\ead{yanchuk@math.tu-berlin.de}

\address{Institute of Mathematics, Technische Universität Berlin, Straße des
17. Juni 136, 10623, Berlin, Germany.}

\author{Т.E. Vadivasova}

\ead{vadivasovate@yandex.ru}

\address{Department of Physics, Saratov State University,83 Astrakhanskaya
Street, Saratov, 410012, Russia}
\begin{abstract}
We present an index for the local sensitivity of spatiotemporal structures
in coupled oscillatory systems based on the asymptotic scaling of
local-in-space, finite-time Lyapunov Exponents. For a system of nonlocally-coupled
R\"{o}ssler oscillators, we show that deviations of this index reflect
the sensitivity to noise and the onset of spatial chaos for the patterns
where coherence and incoherence regions coexist.
\end{abstract}
\maketitle

\section*{Introduction}

The investigation of synchronization and more complex spatiotemporal
structures in coupled oscillatory and spatially extended systems is
a predominant subject of nonlinear science \citep{Pikovsky2001,Kuramoto-book-2002,Afraimovich-1995,Nekorkin-2002,Osipov-2007,Malchow-2007}.
In recent years, after the discovery of structures with coexisting
coherent and incoherent parts, so called chimera states, systems with
nonlocal coupling are of particular interest \citep{Kuramoto-2002,Abrams-2004,Omel2008,AbramsMirolloStrogatzEtAl2008,Laing2010,Omelchenko-2011,Omelchenko2012,Hagerstrom2012,Omelchenko2012a,Martens2013,KapitaniakKuzmaWojewodaEtAl2014,Schmidt2014,Panaggio-2015,Schoell2016,Anishchenko-2016,Semenova-2016,Shepelev-ND-2017,HartBansalMurphyEtAl2016}.
The temporal dynamics of individual elements in such a regime can
be regular (stationary, periodic, quasi-periodic), chaotic, and even
stochastic, see for example \citep{Semenova-2016}. However, despite
of the homogeneous coupling topology, not all elements necessarily
display behavior of the same type. Typically, clusters are formed,
which are characterized by coherent (almost synchronous) behavior
of their constituent elements in time. Elements not belonging to such
clusters exhibit irregular, incoherent behavior.

The chaotic nature of these incoherent ensembles is characterized
by sensitive dependence on initial conditions and is quantitatively
captured by the corresponding local rates of contraction and expansion
along a trajectory, the so called Lyapunov exponents (LEs) \citep{Cesari-1959,Christiansen-1997,Ryabov-2002,Pikovsky2016}.
The maximal LE provides important information about the dynamics of
the system as a whole: chaotic motion, when the maximal LE is positive
and asymptotic convergence to steady state, when it is negative. Modern
computer software facilitates the calculation of the full Lyapunov
spectrum for ensembles consisting of a large number of oscillators
\citep{Popovych-2005,Patra-2016}. For instance, the work \citep{Wolfrum-2011}
studied the full Lyapunov spectrum of a chimera regime, as well as
its dependence on an increasing number of ensemble elements. However,
it is yet unclear how the maximal LE or full Lyapunov spectrum reflect
the dynamical properties of an individual oscillator in a coupled
system. Therefore, special characteristics for the local analysis
in spatially-distributed systems have been considered in the literature.
For instance, local Lyapunov Exponents \citep{Pikovsky-1993} have
been introduced to measure the exponential growth of perturbations
localized in space. However in general, the computation of such local
indicators for all spatially localized perturbations in large ensembles
is numerically challenging.

In this article, we propose the following approach: The finite-time
growth rate of the Lyapunov vector projected onto the subspace corresponding
to a specific oscillator, which we call the \textit{index of local
sensitivity} (ILS). The ILS, in order words, measures the sensitivity
of individual oscillators to external perturbations in a coupled system.
It is computationally cheap in comparison, as it can be computed for
all coupled elements in parallel. This article is organized as follows:
Firstly, we introduce the needed mathematical concepts in more detail.
Secondly, we show how the ILS performs for a system of nonlocally
coupled Rössler oscillators in comparison with the conventional maximal
Lyapunov Exponent. We study the relation between the spatial distribution
of ILSs and the response of different oscillator ensembles to short-time
noise. In particular, we show that elements with larger ILS-values
are more sensitive to such perturbations, as it is expected. In addition,
we reveal that the onset locus of spatial chaos can be characterized
by the ILS: The incoherent part of the chimera state has larger ILS
than the conventional maximum Lyapunov Exponent. 

\section{Index of local sensitivity: Definition and basic properties}

For a given solution $\mathbf{x}$ of the system 
\[
\dot{\mathbf{x}}(t)=f(\mathbf{x}(t)),\quad\mathbf{x}(t)\in\mathbb{R}^{n},
\]
the temporal evolution of a small perturbation vector $\boldsymbol{\xi}(t)$
from $\mathbf{x}(t)$ is characterized by the corresponding Lyapunov
exponent (LE) $\lambda$ \citep{Pikovsky2016}. The LE is defined
as 
\begin{equation}
\lambda:=\limsup\limits _{T\to\infty}\Lambda(t_{0},T),\label{eq:LE}
\end{equation}
where
\begin{equation}
\Lambda(t_{0},T):=\dfrac{1}{T}\ln\dfrac{\|\boldsymbol{\xi}(t_{0}+T)\|}{\|\boldsymbol{\xi}(t_{0})\|},\label{eq:LE-1}
\end{equation}
is the finite-time growth rate of a small perturbation $\boldsymbol{\xi}(t_{0})$
from $\mathbf{x}(t)$ at time $t_{0}$ governed by the linearized
system 
\[
\dot{\boldsymbol{\xi}}(t)=Df(\mathbf{x}(t))\boldsymbol{\xi},\quad\boldsymbol{\xi}(t)\in\mathbb{R}^{n}.
\]
The spectrum of Lyapunov exponents consists of (up to) $n$ numbers
$\lambda_{1}\ge\lambda_{2}\ge,\dots,\ge\lambda_{n}$ corresponding
to different initial perturbations $\boldsymbol{\xi}(t_{0})$, see
e.g. \citep{Wiggins2003} for details. In particular, $\lambda_{1}=\max_{i=1,\dots n}\lambda_{i}$
is called the \emph{maximal Lyapunov exponent,} $\lambda_{\max}$
for short. For a generic initial perturbation $\boldsymbol{\xi}(t_{0})$,
the definition (\ref{eq:LE}) leads to the maximal LE. To remind the
reader, when $\lambda_{\max}$ is positive, the system displays sensitive
dependence on initial conditions and is chaotic. If alternatively,
the maximal LE is negative, any small perturbation asymptotically
converges to $\mathbf{x}$. However, this convergence can be very
slow and preceded by long chaotic transients, as the finite-time growth
rates $\Lambda(t_{0},T)$ can still be positive for large intervals
of time \citep{Politi2010}. In this way, the finite-time growth
rates provide additional information about the dynamics of $\mathbf{x}$.

When we consider an ensemble of oscillators consisting of $N$ elements,
each of them described by $k$ independent variables, the full system
has a phase space of dimension equal to $n=Nk$. Here,  the $Nk$-dimensional
perturbation vector $\boldsymbol{\xi}(t)$ can be represented as $\boldsymbol{\xi}(t)=\left(\boldsymbol{\xi}_{1}(t),\dots,\boldsymbol{\xi}_{N}(t)\right)$,
where $\boldsymbol{\xi}_{i}(t)\in\mathbb{R}^{k},\,i\in\left\{ 1,2,\dots,N\right\} $
is the projection to the subspace corresponding to the $i$-th oscillator.
We propose a new characteristic measuring the response of oscillator
$i$ to a homogeneous perturbation of the form $\boldsymbol{\xi}(t_{0})=\left(\boldsymbol{\xi}_{1}(t_{0}),\dots,\boldsymbol{\xi}_{1}(t_{0})\right)$:
The finite-time growth rate of the Lyapunov vector projected onto
the subspace corresponding to this oscillator, which we call the \textit{index
of local sensitivity} (ILS) and which we define as
\begin{equation}
\Lambda_{i}(t_{0},T):=\dfrac{1}{T}\ln\dfrac{\|\boldsymbol{\xi}_{i}(t_{0}+T)\|}{\|\boldsymbol{\xi}_{i}(t_{0})\|}=\dfrac{1}{T}\ln\sqrt{N}\dfrac{\|\boldsymbol{\xi}_{i}(t_{0}+T)\|}{\|\boldsymbol{\xi}(t_{0})\|}.\label{LIS}
\end{equation}
The ILS is related to the finite-time LE, as 
\[
e^{2\Lambda(t_{0},T)T}=\dfrac{1}{N}\sum_{i=1}^{N}e^{2\Lambda_{i}(t_{0},T)T}
\]
and in some sense, it can be represented as a mean of the individual
contributions of the ILSs to the maximal LE
\begin{equation}
\lambda_{\max}=\limsup_{T\to\infty}\dfrac{1}{2T}\ln\dfrac{1}{N}\sum_{i=1}^{N}e^{2\Lambda_{i}(t_{0},T)T}.\label{LIS4}
\end{equation}
By replacing $\Lambda(t_{0},T)$ by $\lambda_{\max}$ here, we assumed
that the homogeneous perturbation is ''generic'' in the sense that
it is not contained in the subspace corresponding to smaller Lyapunov
exponents $\lambda_{2},\dots,\lambda_{n}$. In the specific case,
when all oscillators are synchronized, all indices of local sensitivity
are identical ($\Lambda_{0}(t_{0},T)=\Lambda_{1}(t_{0},T)=\dots=\Lambda_{N}(t_{0},T)$),
and it holds 
\begin{equation}
\limsup_{T\to\infty}\Lambda_{i}(t_{0},T)=\lambda_{\max}.\label{eq:convergence}
\end{equation}
Even though the ILS can asymptotically converge to the maximal Lyapunov
exponent for less coherent ensembles, the direction (from above or
from below) and the speed provide valuable information about the finite-time
''sensitivity'' of a particular oscillator. In the following sections,
we numerically investigate the dependence of the ILS on the index
$i$, as well as on the reference time $T$. For brevity, from now
on, we omit the dependence on $t_{0}$ and when numerical results
are presented, the associated spatial profile $\mathbf{x}_{i}(t_{0})$
will be clearly stated, where it is appropriate.

\section{Nonlocally coupled Rössler oscillators}

To showcase the approach, we use a system of $N$ nonlocally coupled
chaotic Rössler oscillators, which are known to demonstrate various
spatiotemporal structures, including chimera states \citep{Omelchenko-2012}.
The system under consideration is described by the following set of
$3N$ ordinary differential equations: 
\begin{equation}
\begin{array}{l}
\dot{x}_{i}(t)=-y_{i}(t)-z_{i}(t)+\dfrac{\sigma}{2P}\sum\limits _{k=i-P}^{i+P}\left(x_{k}(t)-x_{i}(t)\right),\\
\dot{y}_{i}(t)=x_{i}(t)+ay_{i}(t)+\dfrac{\sigma}{2P}\sum\limits _{k=i-P}^{i+P}\left(y_{k}(t)-y_{i}(t)\right),\\
\dot{z}_{i}(t)=b+z_{i}(t)(x_{i}(t)-c)+\dfrac{\sigma}{2P}\sum\limits _{k=i-P}^{i+P}\left(z_{k}(t)-z_{i}(t)\right),\\[12pt]
x_{i+N}(t)=x_{i}(t),~~y_{i+N}(t)=y_{i}(t),\\[6pt]
z_{i+N}(t)=z_{i}(t),~~i=1,...N,
\end{array}\label{eq:Rossler}
\end{equation}
such that the underlying coupling topology of the system is periodic,
i.e. it is a nonlocal ring, where index $i$ determines the position
of an oscillator, which is coupled to its $P$-nearest neighbors from
each side with coupling strength $\sigma$. The parameters $a$, $b$,
and $c$ determine the dynamics of an individual oscillator. For numerical
purposes, we confine ourselves to $a=0.2$, $b=0.2$, and $c=4.5$,
so that in the uncoupled case the dynamics of each element is chaotic.
We consider an ensemble consisting of $N=300$ oscillators, each one
coupled to its $P=100$ nearest neighbors. In the following, we numerically
compute the ILS for different values of the coupling strength and
relate the results to the temporal dynamics of the system, in order
to familiarize the reader with this approach.

\subsection{Complete incoherence versus complete chaotic synchronization}

\begin{figure}[!tbh]
\noindent\begin{minipage}[t]{1\linewidth}%
\begin{minipage}[t]{0.45\columnwidth}%
\includegraphics[width=1\linewidth]{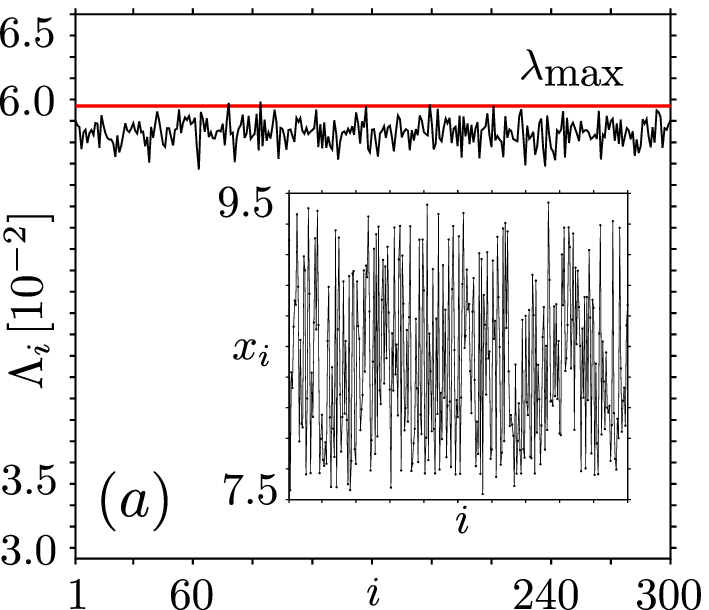}%
\end{minipage}\hfill{}%
\begin{minipage}[t]{0.45\linewidth}%
\includegraphics[width=1\linewidth]{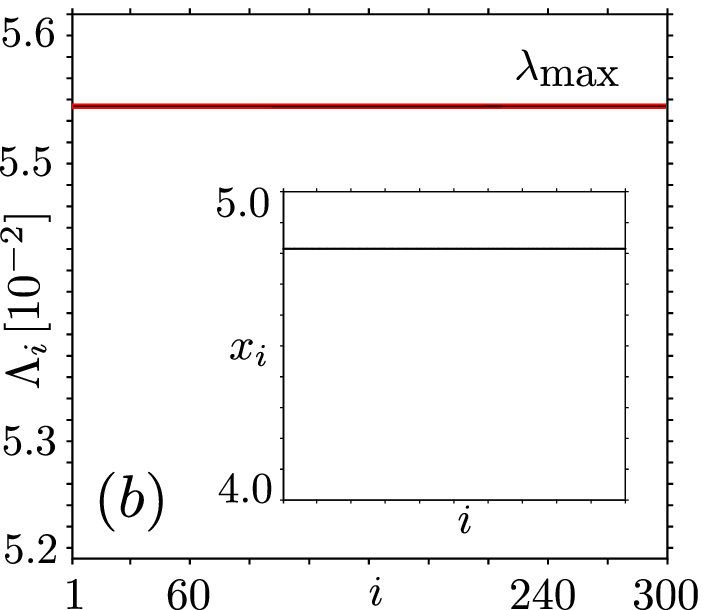}%
\end{minipage}%
\end{minipage}

\caption{\label{pic:uncoupled}(Color online) Spatial distribution of ILS $\Lambda_{i}(T),\,T=10^{4}$
in the regime of complete spatial incoherence (a) $\sigma\geq0$ small,
and complete coherence (b) $\sigma\gg1$. Panels show $\Lambda_{i}(T)$
versus $i$ in comparison with the maximal LE $\lambda_{\text{max}}$
(red). The inset in (a), respectively (b), shows the spatial distribution
of $x_{i}$ at a fixed instant of time. }
\end{figure}

At first we contrast two limiting cases: Complete spatial incoherence
and the regime of complete synchronization. For small positive values
of the coupling strength $\sigma\approx0$, one observes complete
spatial incoherence \citep{Politi2006,Rosenblum1996}. Here for randomly
chosen initial conditions, each individual element is chaotic,  and
 at a fixed moment in time the spatial distribution does not exhibit
any apparent coherent ensemble. The corresponding ILSs $\Lambda_{i}(T)$
are shown in Fig.~\eqref{pic:uncoupled}(a) alongside an example
of spatial distribution of the oscillators $x_{i}$, see Fig.~\eqref{pic:uncoupled}(a,
inset).
 Naturally in this case, the ILS are close to the maximal LE with
some local variation among the oscillators, as finite-time LE do not
depend continuously on the initial condition. Note that in this case
all $\Lambda_{i}$ attain smaller values than the maximal LE, as it
corresponds to the direction along which the expansion is maximal.
Increasing the coupling strength $\sigma$ to large values, complete
chaotic synchronization can be observed in the system \citep{Pecora1998,Pecora1997,Yanchuk2001a,Yanchuk2003}.
In this regime, after some transient, all elements still oscillate
chaotically, but synchronously so, with $\mathbf{x}_{i}(t)=\mathbf{x}_{j}(t)$,
for all $i,j$. As the dynamics of the oscillators is identical, all
$\Lambda_{i}(T)$ coincide and as $T\to\infty$, they converge to
the maximal Lyapunov exponent up to the precision given by the resolution
of Fig.~\ref{pic:uncoupled}(b). For a corresponding example spatial
distribution of the oscillators $x_{i}$ see Fig.~\eqref{pic:uncoupled}(b,
inset). We have made a simple observation here: incoherent (coherent)
ensembles of oscillators posses incoherent (coherent) ILS. Throughout
this paper, we will exploit this fact and investigate how the value
and scaling of ILSs reflect dynamical properties of the system state
and what additional information can be obtained from the ILS. 

\subsection{Partial chaotic synchronization\label{subsec:partial}}

\begin{figure}[!tbh]
\includegraphics[width=1\linewidth]{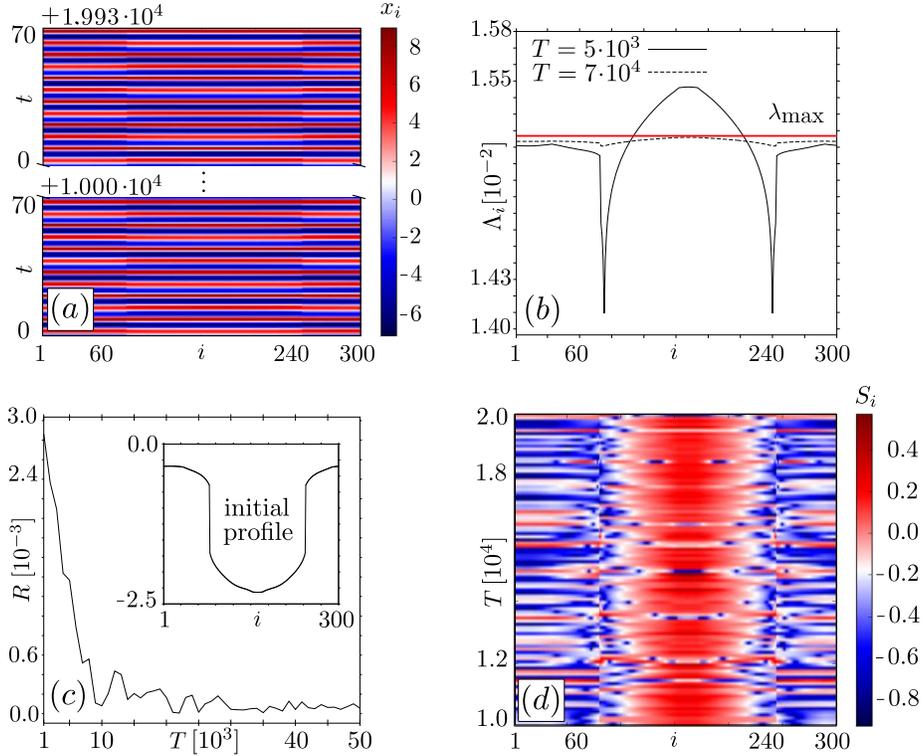}\caption{\label{pic:part_coh}(Color online) Spatial distribution of $\Lambda_{i}(T)$
in the regime of partial synchronization of chaotic dynamics for intermediate
coupling strength $\sigma=0.05$. Panel (a) shows an example trajectory
of the system for $t\in[10^{4},2\cdot10^{4}]$ in a spatiotemporal
plot. Panel (b) shows the corresponding spatial distribution of the
ILS $\Lambda_{i}(T)$ for $T=5\cdot10^{3}$ (solid line), $T=7\cdot10^{4}$
(dotted line) and the maximal LE $\lambda_{\text{max}}$ of the full
system along this solution. One observes that $\Lambda_{i}(T)$ converges
to $\lambda_{\text{max}}$ and the local minima along the spatial
distribution correspond to the boundaries between the two clusters.
Panel (c) shows the spatial range of ILSs $R_{\text{ILS}}(T)=\max_{i}\Lambda_{i}(T)-\min_{i}\Lambda_{i}(T)$
versus the reference time T for a fixed initial profile (inset). Panel
(d) shows the temporal evolution of the rescaled ILS $S_{i}(T)=(\Lambda_{i}(T)-\Lambda(T))/R_{\text{ILS}}(T)$.}
\end{figure}

For intermediate values of the coupling, one can observe partial coherence
between the individual chaotic oscillators \citep{Omelchenko-2012}.
Solutions in this regime are characterized by piecewise smooth instantaneous
spatial profiles, see Fig.~\ref{pic:part_coh}(a), where almost all
adjacent elements of the ensemble oscillate approximately synchronously,
but two distant elements can have (possibly very) different instantaneous
states. For example, the solution shown in Fig.~\ref{pic:part_coh}(a)
is almost coherent in space: small (large) values of oscillator $\mathbf{x}_{i}$
correspond to small (large) values of $\mathbf{x}_{j}$ for all $j$.
However, it exhibits two points of discontinuity in its spatial profile
at the oscillators $i=85$, and $i=240$ respectively, such that it
can be thought of as two clusters of oscillators $i=1,\dots,85,240,\dots,300$
and $i=85,\dots,240$ each one coherent in space. 

The distribution of ILSs $\Lambda_{i}(T)$ along this solution is
shown in Fig.~\ref{pic:part_coh}(b). We observe that it varies continuously
 in space and has a pronounced maximum about oscillator $i=160$,
where $\Lambda_{i}(T)>\lambda_{\text{max}}$, which corresponds to
local rates of expansion higher than the maximum Lyapunov exponent.
This observation contrasts with our intuition that the maximal sensitivity
of the ensemble must be observed in the regions around the profile
distortions. Surprisingly, the minima of the ILS distribution are
observed around the points, where the spatial profile is discontinuous. 

Our numerical simulations show that the range of $\Lambda_{i}(T)$,
that is $R_{\text{ILS}}(T):=\max_{i}\Lambda_{i}(T)-\min_{i}\Lambda_{i}(T)$
decreases with time \ref{pic:part_coh}(c) and asymptotically, all
sensitivity indices $\Lambda_{i}(T)$ have approximately the same
value $\lambda_{\max}$ Fig.~\ref{pic:part_coh}(b). The qualitative
explanation for this effect is as follows: for $T\rightarrow\infty$
the vector of perturbation of the whole system line up with the direction
corresponding to the maximal Lyapunov exponent \citep{Pikovsky2016,GinelliPoggiTurchiEtAl2007}.
This perturbation is changing along the phase trajectory, however,
the lengths of all the projections of a perturbation vector in partial
oscillators are changed proportionally. In the limit, any projection
either increases or decreases exponentially with the same rate, which
is determined by the value of the maximal Lyapunov exponent $\lambda_{\max}$. 

We emphasize that the variation of ILS for different oscillators is
an effect of the finite calculation time $T$. In order to use this
characteristic, one should choose it in an optimal way: small enough
to avoid the asymptotic limit (in our case study of the Rössler system
$T<10^{4}$), and large enough to exceed the characteristic timescale
of the system. Interestingly, we observe that, although the local
variation in $\Lambda_{i}(T)$ decreases, the local maxima and minima
of the spatial ILS profile remain over time, see Fig.~\ref{pic:part_coh}(d).
The ILS distribution rescaled by $R_{\text{ILS}}$ forms a characteristic
shape providing qualitative and quantitative information about the
sensitivity of different oscillators in the ensemble.

\section{Index of local sensitivity: Results}

\subsection{Local sensitivity to noise}

\begin{figure}[!t]
\noindent\begin{minipage}[t]{1\linewidth}%
\includegraphics[width=1\linewidth]{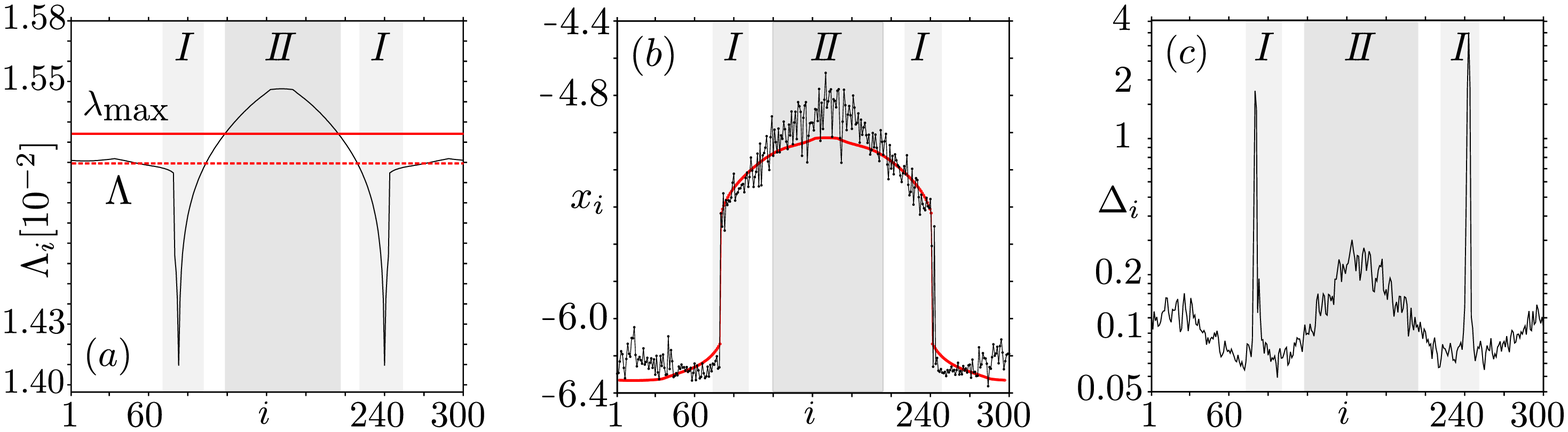}%
\end{minipage}\caption{\label{pic:deviation}(Color online) Influence of uniform noise with
low intensity $D=10^{-5}$ on the spatial structure in the regime
of partial synchronisation $\sigma=0.05$, see also Fig.~\ref{pic:part_coh}.
Panel (a) shows ILS distribution $\Lambda_{i}(T)$ versus $i$ for
$T=5\cdot10^{3}$ (black, solid), the finite-time LE $\Lambda_{i}(T)$
for $T=5\cdot10^{3}$ (red, dashed) and the maximal LE $\lambda_{\text{max}}$
(red, solid). Spatial Regions I (Oscillators $i$ with $\Lambda_{i}(T)<\Lambda(T)$)
and II (Oscillators $i$ with $\Lambda_{i}(T)>\lambda_{\text{max}}$)
are highlighted in every panel. Panel (b) shows the influence of noise
(black) on a spatial profile (red). Panel (c) shows the degree of
incoherence $\Delta_{i}$ averaged over $t\in[5\cdot10^{3},7\cdot10^{4}]$.}
\end{figure}

In this section, we study the impact of a  noisy perturbation  on
a partially coherent solution as shown in Fig.~\eqref{pic:part_coh}(a).
We relate our findings to the ILS distribution (Fig.~\eqref{pic:deviation}(a))
and reveal how the deviations in $\Lambda_{i}(T)$ influence the response
of an individual oscillators. In particular, we show that the oscillators
$i$ with small ILS correspond to stable (with respect to noise) regions
of the profile, see Fig.~\eqref{pic:deviation}(b, region I). On
the other hand, oscillators with ILS locally higher than the maximal
LE, $\Lambda_{i}(T)>\lambda_{\text{max}}$, respond to external perturbation
more strongly, see Fig.~\eqref{pic:deviation}(b, region II). 

To do so, we consider the following stochastic system 
\begin{equation}
\begin{array}{l}
\dot{x}_{i}=-y_{i}-z_{i}+\dfrac{\sigma}{2P}\sum\limits _{k=i-P}^{i+P}\left(x_{k}-x_{i}\right)+\sqrt{2D(i,t)}n(t),\\
\dot{y}_{i}=x_{i}+ay_{i}+\dfrac{\sigma}{2P}\sum\limits _{k=i-P}^{i+P}\left(y_{k}-y_{i}\right)+\sqrt{2D(i,t)}n(t),\\
\dot{z}_{i}=b+z_{i}(x_{i}-c)+\dfrac{\sigma}{2P}\sum\limits _{k=i-P}^{i+P}\left(z_{k}-z_{i}\right)+\sqrt{2D(i,t)}n(t),\\[12pt]
\end{array}\label{eq:Rossler_noise}
\end{equation}
with periodic boundary conditions, as in (\ref{eq:Rossler}). Here,
$n(t)$ is a normalized source of Gaussian white noise with intensity
$D(i,t)$. 

In order to measure the effect of noise applied to the system, we
use the  quantity
\begin{equation}
\Delta_{i}=\langle(2x_{i}(t)-x_{i+1}(t)-x_{i-1}(t))^{2}\rangle\label{eq:deviation}
\end{equation}
as the degree of the local incoherence of the spatial profile at point
$i$, see \citep{Shepelev-ND-2017,Shepelev-2017-FHN} for more details.
Here, the averaging $\langle\cdot\rangle$ is performed with respect
to time. This characteristic displays the averaged ''curvature'',
as the expression $\left(2x_{i}(t)-x_{i+1}(t)-x_{i-1}(t)\right)^{2}$
is a measure of the local deviation from the linear state. $\Delta_{i}$
admits small values for coherent states and larger values for incoherent
as it is shown in \citep{Shepelev-ND-2017}. 

At first, we consider a constant noise intensity $D(i,t)\equiv D$,
uniform for all time $t$ and oscillators $i$. Figure~\ref{pic:deviation}(b)
shows an example profile in space for a fixed instant of time with
and without noise. The corresponding distribution of $\Delta_{i}$
in space is shown in Fig.~\ref{pic:deviation}(c), where the time
averaging is performed for $5\cdot10^{3}$ time units. The distinct
peaks in Fig.~\ref{pic:deviation}(c, region I) correspond to the
profile distortion. At the same time, the maximum of $\Delta_{i}$
in region II, coincides with the maximum in the ILS distribution in
the plots in Figs.~\ref{pic:deviation}(a). Thus, the elements with
the largest values of $\Lambda_{i}(T)$ are most sensitive to the
noise influence, in contrast to the elements, which are characterized
by smaller values of ILS.

\begin{figure}[!ht]
\noindent\begin{minipage}[t]{1\linewidth}%
\includegraphics[width=1\linewidth]{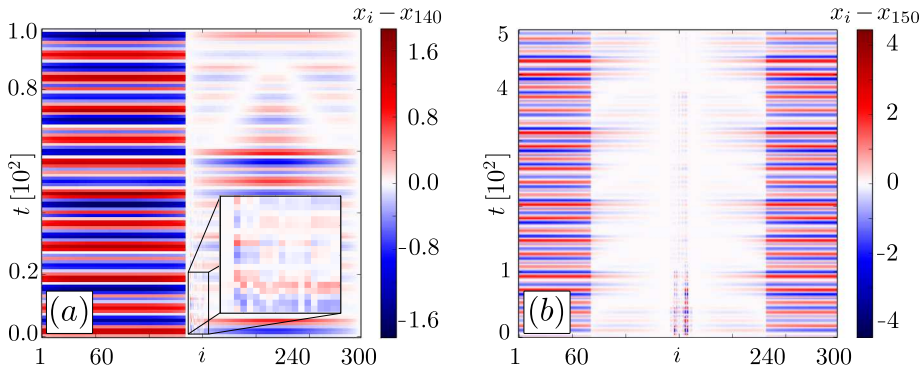}%
\end{minipage}

\caption{\label{pic:loc_noise_snap}(Color online) Influence of spatially localized short-time noise with intensity $D=0.05$
in the regime of partial synchronization of chaotic dynamics for intermediate coupling strength $\sigma=0.05$.
Panels (a) and (b) show the spatiotemporal evolution of a perturbation in regions I (a) and II (b) (see Fig. 3 for the definition of regions I and II) in a corresponding reference frame: $x_{i}-x_{140}$
  in (a) and $x_{i}-x_{150}$  in (b). In panel (a), oscillator numbers are shifted such that the perturbation appears centered. }
\end{figure}

Additionally, when studying the influence of localized short-time
noise, 
\begin{equation}
\begin{array}{l}
D(i,t)=\begin{cases}
D,~~i\in[i_{1};i_{2}]\,\,\mbox{and}\,\,t\in[0,T_{n}],\\
0,~~\,\,\mbox{otherwise,}
\end{cases}\end{array}\label{eq:ext_noise}
\end{equation}
where $\left[i_{1},i_{2}\right]$ is the spatial and $[0,T_{n}]$
the temporal interval of the noise action, we can relate the value
of the ILS, to the characteristic decay time of this perturbation.
We consider the system (\ref{eq:Rossler_noise}) without noise to
be in the regime of partial coherence, as studied above (Fig.~\ref{pic:part_coh}),
with noise parameters chosen as $D=0.05$ and $T_{n}=0.1$. Figure~\eqref{pic:loc_noise_snap}
exhibits the spatiotemporal plot of two local perturbations applied
in regime I (low ILS, Fig.~\eqref{pic:loc_noise_snap}(a)), and regime
II (high ILS, Fig.~\eqref{pic:loc_noise_snap}(b)) respectively.
One can observe that the perturbation in the high-sensitive region
II persists significantly longer then in the low-sensitive region
I. Thus, the ILS can serve as a sensitivity measure to noisy perturbations.
The non-coherent response, which is induced by the localized short-time
high-intensity influence of noise, persists significantly longer in
the case when the perturbation is applied to a high-sensitive region
(with higher ILS) than in case of the influence to the region with
low values of ILS. The obtained results indicate the non-homogeneity
of the response of system (\ref{eq:Rossler}) in the regime of partial
coherence.

\subsection{Local sensitivity of chimera states}

In the regime of partial chaotic synchronization, so called ''chimera''
states have been observed. There are 2 main flavors in which a spatially
coherent profile can become locally incoherent \citep{Anishchenko-2016}:
(1) An oscillator that is close to the boundary of two coherent clusters,
say cluster 1 and 2, can be irregularly ``shifted'' in space with
respect to the instantaneous phase of adjacent oscillators belonging
to cluster 1, but the oscillator itself has values comparable to an
oscillator in cluster 2: the so-called ``phase'' chimera, which
needs to be distinguished from a ``chimera of phase oscillators''.
(2) In an ``amplitude'' chimera, the  oscillators locally possess
an irregular spatial distribution with respect to the surrounding
cluster. So naturally, phase chimeras occur in the region of low ILS
(compare region I in Fig.~\ref{pic:deviation}(a)) and amplitude
chimeras in the region of high ILS (compare region II in Fig.~\ref{pic:deviation}(a)).
In this section, we numerically investigate the distributions of ILS
for these regimes. 

\begin{figure}[!ht]
\noindent\begin{minipage}[t]{1\linewidth}%
\includegraphics[width=1\linewidth]{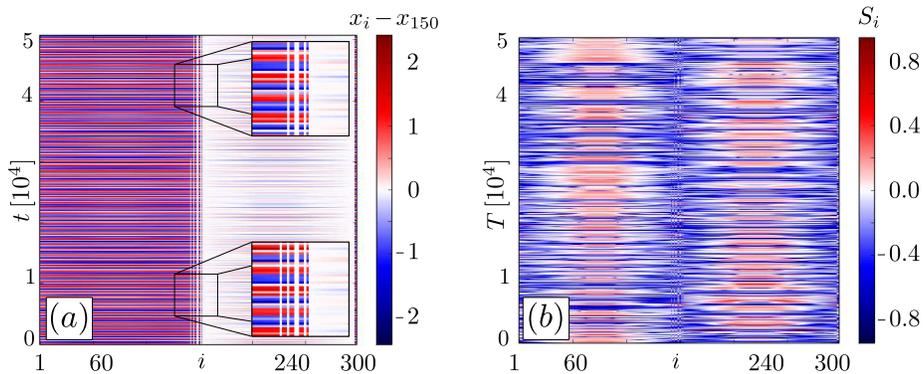}%
\end{minipage}

\caption{\label{pic:phase_chimera}(Color online) Spatial distribution of $\Lambda_{i}(T)$
in the regime of partial coherence (``phase'' chimera) of chaotic
dynamics for coupling strength $\sigma=0.044$. Panel (a) shows the
temporal evolution of the oscillators $x_{i}$ with respect to oscillator
$150$ that is part of the incoherent ensemble. Panel (b) shows the
temporal evolution of the rescaled ILS $S_{i}(T)=(\Lambda_{i}(T)-\Lambda(T))/R_{\text{ILS}}(T)$.
(In panels (a) and (b) oscillator numbers have been shifted, such
that the chimera state appears centered)}
\end{figure}

Firstly we study the phase chimera shown in Fig.~\ref{pic:phase_chimera}(a).
Here, elements of the incoherent cluster oscillate almost periodically
in time. Their instantaneous amplitudes are almost identical (changes
smoothly with $i$), and the adjacent elements oscillate either in-phase
or anti-phase. Moreover, these shifts are irregularly distributed
in space. It has been argued that this chimera type is stable in time
and resistant to external influence \citep{Semenova-2017}, which
is also reflected by the ILS. Figure~\ref{pic:phase_chimera}(b)
shows the distribution of ILSs $\Lambda_{i}(i)$ rescaled by their
range $R_{\text{ILS}}(T)$ corresponding to the phase chimera in (a).
The ILS distribution in the coherence cluster changes smoothly leading
to the similar sensitivity of adjacent elements. On the contrary,
the ILS of the incoherence cluster varies irregularly and most importantly
has lower ILS than the oscillators in the coherent clusters. 

\begin{figure}[!ht]
\noindent\begin{minipage}[t]{1\linewidth}%
\includegraphics[width=1\linewidth]{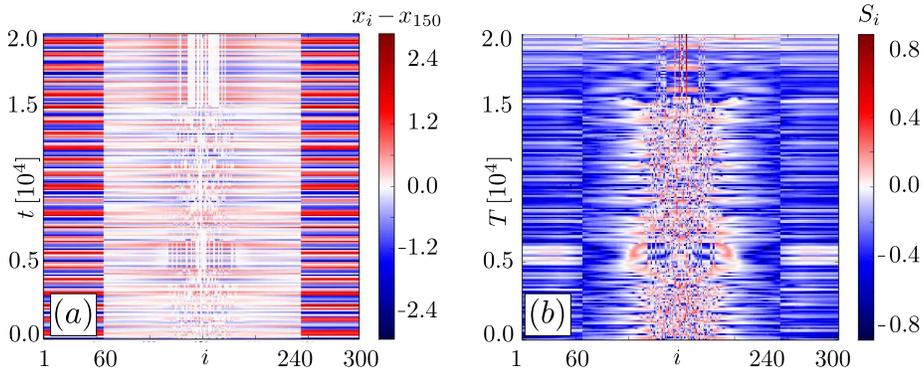}%
\end{minipage}\caption{\label{pic:ampl_chimera}(Color online) Spatial distribution of $\Lambda_{i}(T)$
in the regime of partial coherence (``amplitude'' chimera) of chaotic
dynamics for coupling strength $\sigma=0.04$. Panel (a) shows the
temporal evolution of the oscillators $x_{i}$ with respect to oscillator
$150$ that is part of the incoherent ensemble. Panel (d) shows the
temporal evolution of the rescaled ILS $S_{i}(T)=(\Lambda_{i}(T)-\Lambda(T))/R_{\text{ILS}}(T)$.
(In panels (a) and (b) oscillator numbers have been shifted, such
that the chimera state appears centered)}
\end{figure}

An example of amplitude chimera is presented in Fig.~\ref{pic:ampl_chimera}(a).
Amplitude chimeras in ensembles of chaotic oscillators are presumably
metastable states, unlike for the phase chimeras. They exist for a
finite (but possibly very long) time \citep{Semenova-2017} and are
regarded as sensitive to perturbations in the form of noise.

Fig.~\ref{pic:ampl_chimera}(b) represents the spatial distribution
of $\Lambda_{i}(T)$ rescaled by range $R_{\text{ILS}}(T)$ in the
regime of the amplitude chimera. The ILS distribution is continuous
in the coherence cluster, while the values of $\Lambda_{i}(T)$ can
significantly vary for the elements of the incoherence cluster. The
ILS has a maximum in the region of an incoherence cluster, unlike
the phase chimera. Thus, the region of an amplitude chimera is the
least stable part of the spatial structure. 

\section*{Conclusion}

We are convinced that the index of local sensitivity enriches the
numerical toolbox for the study of spatiotemporal structures and the
fast growing field of network analysis. There are two major reasons
for this: 

(1) It has a clear interpretation in terms of the finite-time growth
rate of perturbations to a specific oscillator. We have shown that
this interpretation is indeed valid, as the value of ILS is related
to the decay time of perturbations by short-time noise applied to
a specific oscillator. In this way, we are able to identify elements,
which are most sensitive to small perturbations (including noise).
This information is important for the development and evaluation of
methods to influence single oscillators and larger ensembles by an
external control and to improve the structure of a general network
in order to render it more robust to external perturbations. From
a theoretical point of view, we showed that the ILS provides valueable
information in the study of complex spatial structures, as it can
predict the onset locus of spatial chaos in the case of chimera states
in a system of coupled R\"{o}ssler oscillators. 

(2) It is computationally cheap, such that it can be applied to coupled
systems with a large number of elements. As a finite-time measurement,
it can be calculated for a varying range of timescales that, needless
to say, need to be adjusted to the phenomenon under consideration.
We plan to investigate the scaling behavior of the ILS as $T\to\infty$
in more detail, as for the solutions we considered, each ILS converges
to the maximum Lyapunov Exponent of the full system (up to numerical
intergration error). Here, many questions arise: Does the distribution
of ILSs rescaled by $R_{ILS}$ converge to an asymptiotic shape? Does
negativity of an ILS imply local contraction and the spatiotemporal
structure becomes fixed in time? For a general system, this of course
need not necessarily be the case and there are many such open questions
 that can be investigated in the future. 

\section*{Acknowledgments}

This work was supported by the German Research Foundation (DFG) in
the framework of the Collaborative Research Center (SFB) 910, projects
A3 and B11.  T. E. V. acknowledges support from the Russian Science
Foundation (grant No.~16-12-10175).


\bibliography{sy-bibliography}
\bibliographystyle{unsrt}

\end{document}